\documentclass[9pt]{article}
\usepackage{latexsym, amsfonts}
\newcommand{\be}{\begin{equation}}
\newcommand{\ee}{\end{equation}}
\newcommand{\bea}{\begin{eqnarray}}
\newcommand{\eea}{\end{eqnarray}}
\newcommand{\bean}{\begin{eqnarray*}}
\newcommand{\eean}{\end{eqnarray*}}
\newcommand{\brray}{\begin{array}}
\newcommand{\erray}{\end{array}}
\newcommand{\ben}{\begin{equation}{nonumber}}
\newcommand{\een}{\end{equation}{nonumber}}

\newtheorem{dfn}{Definition}[section]
\newtheorem{thm}[dfn]{Theorem}
\newtheorem{lmma}[dfn]{Lemma}
\newtheorem{ppsn}[dfn]{Proposition}
\newtheorem{crlre}[dfn]{Corollary}
\newtheorem{xmpl}[dfn]{Example}
\newtheorem{rmrk}[dfn]{Remark}
\newcommand{\bdfn}{\begin{dfn}}
\newcommand{\bthm}{\begin{thm}}
\newcommand{\blmma}{\begin{lmma}}
\newcommand{\bppsn}{\begin{ppsn}}
\newcommand{\bcrlre}{\begin{crlre}}
\newcommand{\bxmpl}{\begin{xmpl}}
\newcommand{\brmrk}{\begin{rmrk}}
\newcommand{\edfn}{\end{dfn}}
\newcommand{\ethm}{\end{thm}}
\newcommand{\elmma}{\end{lmma}}
\newcommand{\eppsn}{\end{ppsn}}
\newcommand{\ecrlre}{\end{crlre}}
\newcommand{\exmpl}{\end{xmpl}}
\newcommand{\ermrk}{\end{rmrk}}

\newcommand{\IT}{\mathbb{T}}

\newcommand{\IZ}{\mathbb{Z}}

\newcommand{\cla}{{\cal A}}

\newcommand{\clg}{{\cal G}}
\newcommand{\clh}{{\cal H}}

\newcommand{\cll}{{\cal L}}

\newcommand{\clq}{{\cal Q}}

\def\a*{{\cal A}_{h,*}}
\def\B{{\cal B}(h)}
\def\B1{{\cal B}_1(h)}
\def\b{{\cal B}^{\rm s.a.}(h)}
\def\b1{{\cal B}^{\rm s.a.}_1(h)}

\newcommand{\ot}{\otimes}

\newcommand{\raro}{\rightarrow}

\def \qed {$\Box$}

\begin{document}
	\[
\]
\begin{center}
{\large {\bf Quantum Isometry Group of the n tori}}\\
by\\
{\large Jyotishman Bhowmick} {\footnote {The support from National Board of Higher Mathematics,  India,
 is gratefully acknowledged.}} \\
{\large Stat-Math Unit, Kolkata Centre,}\\
{\large Indian Statistical Institute}\\
{\large 203, B. T. Road, Kolkata 700 108, India}\\
{e mail: jyotish\_\ r@isical.ac.in }\\
\end{center}
\begin{abstract}
    
    We show that the Quantum Isometry Group (  introduced in \cite{goswami} ) of the n tori $ \IT^{n} $ coincides with its classical isometry group i.e, there does not exist any faithful isometric  action on $ \IT^{n} $ by a genuine ( noncommutative as a $ C^{\ast} $ algebra ) compact quantum group . Moreover, using a result in \cite{bhowmick goswami}, we conclude that the Quantum Isometry group of the noncommutative n tori is a Rieffel deformation of the Quantum Isometry Group of the commutative n tori.
      
  \end{abstract} 
  \section{Introduction}
 In \cite{goswami}, Goswami has defined the quantum isometry group of a noncommutative manifold (given by spectral triple), motivated by the definition and study of 
quantum permutation groups of finite sets and finite graphs by a number of mathematicians (see, e.g.
\cite{ban1}, \cite{ban2}, \cite{wang}, and references therein) and using some ideas of Woronowicz and Soltan ( see \cite{soltan} ).

      The main ingredient of this theory is the Laplacian $\cll$ coming from the spectral triple $ (\cla^\infty, \clh, D) $ satisfying certain regularity conditions (see \cite{goswami} for its construction), which coincides with  the Hodge Laplacian $-d^\ast d$ (restricted on space of smooth functions) in the classical case, where $d$ denotes the de-Rham differential.
     
The linear span of eigenvectors of $\cll$, which is a subspace of $\cla^\infty$, is denoted by $\cla^\infty_0$, and it is assumed that $\cla^\infty_0$ is norm-dense in the $C^*$-algebra  $\cla$ obtained by completing $\cla^\infty$.  The $\ast$-subalgebra  of $\cla^\infty$ generated by $\cla^\infty_0$ is denoted by $\cla_0$.Then  $\cll(\cla^\infty_0) \subseteq \cla^\infty_0$,   and a compact quantum group $(\clg,\Delta)$ which has an action $\alpha$ on $\cla$ is said to act smoothly and isometrically on the noncommutative manifold $(\cla^\infty, \clh, D)$ if for every state $ \phi $ on $ \clg $, $({\rm id} \ot \phi) \circ \alpha(\cla^\infty_0) \subseteq \cla^\infty_0$ for all state $\phi$ on $\clg$, and also $({\rm id} \ot \phi) \circ \alpha$ commutes with $\cll$ on $\cla^\infty_{0}$.

 One can then consider the category of all compact quantum groups acting smoothly and isometrically on $\cla$, where the morphisms are quantum group morphisms which intertwin the actions on $\cla$. It is proved in \cite{goswami}(under some regularity assumptions, which are valid for any compact connected Riemannian spin manifold with the usual Dirac operator) that there exists a universal object in this category, and this universal object is defined to be the quantum isometry group of $(\cla^\infty,\clh,D)$, denoted by $QISO(\cla^\infty, \clh, D)$, or simply as $QISO(\cla^\infty)$ or even $QISO(\cla)$ if the spectral triple is understood.
     
 It is important to explicitly describe quantum isometry groups of sufficiently many classical and noncommutative manifolds.In \cite{bhowmick goswami}, the Quantum Isomtery Group of the two torus has been explicitely computed.
 
 Although there are similarities between the basic principle of the computation of \cite{bhowmick goswami} that of  the present article, which is to exploit the eigenspaces of the Laplacian, there is a key difference between the two proofs which is as follows:
 
The computation of the quantum isometry group of  $ \cla_{\theta} $ in Section 2.3 of \cite{bhowmick goswami} is valid for all $ \theta ,$ i.e. all $\theta$ has been 
treated on
 the same footing. On the other hand,  we have adopted a different strategy in the present article based on the observation that
 the proof of \cite{bhowmick goswami} can be simplified quite a bit for the special case of $\theta=0$, i.e. 
the commutative case, and this simplification even goes through easily for any $n \geq 3$.
Thus, in this article we first compute the quantum isometry group of the commutative $ C^* $ algebra $ C ( \IT^n ).$ 
This commutativity assumption decreases a lot of computations because it suffices to show
 that the quantum isometry group is a commutative $ C^* $ algebra and hence has to coincide with the classical isometry group. 
Then  the structure of the quantum isometry group of $ \IT^{n}_{\theta}$ for an arbitrary value of $\theta$ follows immediately from Theorem 3.13 of \cite{bhowmick goswami}.

 Throughout the paper, we have denoted by $\cla_1 \ot \cla_2$ the minimal (injective) $C^*$-tensor product between two $C^*$-algebras $\cla_1$ and $\cla_2$. The symbol $\ot_{\rm alg}$ has been used to denote the algebraic tensor product between vector spaces or algebras.  
 
 For a compact quantum group $\clg$, let the dense unital $\ast$-subalgebra generated by the matrix coefficients of irreducible  unitary representations be denoted by $\clg_0$. The coproduct of $\clg$, say $\Delta$, maps $\clg_0$ into the algebraic tensor product $\clg_0 \ot_{\rm alg} \clg_0$, and there exist canonical antipode and counit defined on $\clg_0$ which make it into a Hopf $\ast$-algebra ( see \cite{woro} for the details ).

\brmrk

\label{T_n_Van_Daele}
      
 In \cite{goswami}, it was assumed that the compact quantum groups are separable. But in \cite{VanDaele2}, the separability assumption for the $ C^{*} $ algebra of the underlying compact quantum group was removed. It can be easily seen that the separability was not at all used in the proofs of \cite{goswami} and hence all the results in \cite{goswami} goes verbatim in the non separable case.      
       
\ermrk

   \subsection{ Quantum isometry group of the commutative  n-tori} 
   
   Consider  $ C ( \IT^n ) $ as the universal commutative $ C^{*} $ algebra generated by n commuting unitaries $ U_{1}, U_{2},.....U_{n}  $.It is clear that  the set $ \{ U_{i}^{m}U_{j}^{n} : m,n \in \IZ \}  $ is  an orthonormal basis for  $ L^{2} ( C ( \IT^n ) , \tau ),     $ where $ \tau$ denotes the unique faithful normalized trace on $ C ( \IT^n ) $ given by, $\tau ( \sum a_{m n} U_{i}^{m} U_{j}^{n} ) = a_{0 0} $ which is just the integration against the Haar measure.     
      We shall denote by  $ \left\langle A , B \right\rangle = \tau ( A^{*} B )$    the inner product on $\clh_0:=L^2(C ( \IT^n ),\tau)$. Let $ {C ( \IT^n )}^{\rm fin}$ be the unital  $\ast$-subalgebra generated by finite complex linear combinations of $U^mV^n$, $m,n \in \IZ$.The Laplacian $\cll$ is given by $\cll(U_{1}^{m_{1}} ......U_{n}^{m_{n}})=-(m_{1}^2 + ...m_{n}^2) U_{1}^{m_{1}}......U_{n}^{m_{n}},$ and it is also easy to see that the algebraic span of eigenvectors of $\cll$ is nothing but the space ${C ( \IT^n )}^{\rm fin}$, and moreover, all the assumptions in \cite{goswami} required for defining the quantum isometry group are satisfied. 
         
         Let $\clq$ be the quantum isometry group coming from the above laplacian, with the 
  smooth isometric action of $\clq $  on $ C ( \IT^n ) $ given by $\alpha : C ( \IT^n ) \raro C ( \IT^n ) \ot \clq$. 
     By definition,  $ \alpha $  must keep invariant the eigenspace of $ \cll $ corresponding to the eigen value $- 1 $ , spanned by $ U_{1},.....U_{n},U_{1}^{-1},.......,U_{n}^{-1} $.Thus, the action $ \alpha $ is given by:   
    $$  \alpha ( U_{i} ) = \sum_{j = 1}^{n} U_{j} \otimes A_{ij} + \sum_{j=1}^{n} U_{j}^{-1} \otimes B_{ij} ,$$
  
   where $ A_{ij},B_{ij} \in \clq ,i,j =1,2....n $.By faithfulness of the action of quantum isometry group (see \cite{goswami}), the norm-closure of the unital $\ast$-algebra generated by $A_{ij},B_{ij} ;i,j=1,2,....n $ must be the whole of $\clq$.
   
   Next we derive a number of conditions on $ A_{ij},B_{ij}, i,j = 1,2,...n $ using the fact that $ \alpha $ is a $ \ast $   homomorphism.
  \blmma    
    \label{Lemma 1}

   The condition $ U^{*} U = 1 = U U^{*} $ gives:
   
  \be \label{lem1.1a}  \sum_{j} ( A^{*}_{ij} A_{ij} +  B^{*}_{ij} B _{ij} ) = 1    \ee

  \be \label{lem1.2} B^{*}_{ij} A_{ik} +  B^{*}_{ik} A_{ij} = 0 ~ \forall j \neq k \ee

  \be \label{lem1.3} A^{*}_{ij} B_{ik} +  A^{*}_{ik} B_{ij} = 0 ~ \forall j \neq k \ee
   
  \be \label{lem1.4} A^{*}_{ij} B_{ij} = B^{*}_{ij} A_{ij} = 0 \ee

  \be \label{lem1.5} \sum_{j} ( A_{ij} A^{*}_{ij} + B_{ij} B^{*}_{ij} ) = 1 \ee
   
   \be \label{lem1.6}A_{ik} B^{*}_{ij} + A_{ij}B^{*}_{ik} = 0 ~ \forall j \neq k  \ee

  \be \label{lem1.7}B_{ik} A^{*}_{ij} + B_{ij}A^{*}_{ik} =  0 ~ \forall j \neq k \ee
   
  \be \label{lem1.8} A_{ij} B^{*}_{ij} = B_{ij}A^{*}_{ij} =  0 \ee

   \elmma
 {\it Proof :}\\
  We get ( \ref{lem1.1a} ) - ( \ref{lem1.4} ) by using the condition $ U_{i}^{*} U_{i} = 1 $ along with the fact that $ \alpha $ is a homomorphism and then comparing the coefficients of $ 1, U_{j}U_{k}, U_{j}^{-1}U_{k}^{-1},$ ( for $ j \neq k ),$ $ U_{j}^{-2}, U_{k}^{-2} . $
   
   Similarly the condition $ U_{i} U_{i}^{*} = 1 $ gives ( \ref{lem1.5} )  -( \ref{lem1.8} ).\qed \vspace{4mm}

  Now, $ \forall i \neq j, U_{i}^{*}U_{j}, U_{i}U_{j}^{*} $ and $ U_{i}U_{j} $ belong to the eigenspace of the laplacian with eigenvalue $-2$, while $ U_{k}^{2}, U_{k}^{-2} $ belong to the eigenspace corresponding to the eigenvalue $-4$.As $ \alpha $ preserves the eigenspaces of the Laplacian, the coefficients of $ U_{k}^{2}, U_{k}^{-2} $ are zero $ \forall k $ in $ \alpha( U_{i}^{*}U_{j} ),\alpha( U_{i}U_{j}^{*} ), \alpha( U_{i}U_{j} ) ~ \forall i \neq j. $

   We use this observation in the next lemma.
   
   \blmma
   \label{Lemma 3}
   
   $ \forall k $ and $ \forall i \neq j ,$
   
   \be \label{lem3.1}  B^{*}_{ik} A_{jk} = A^{*}_{ik} B_{jk} = 0  \ee

   \be \label{lem3.2}  A_{ik}B_{jk} = B_{ik}A^{*}_{jk} = 0  \ee

   \be \label{lem3.3}  A_{ik} A_{jk} = B_{ik} B_{jk} = 0  \ee

   \elmma
   
    {\it Proof :}\\   
   The equation ( \ref{lem3.1} ) is  obtained from the coefficients of $  U_{k}^{2} $ and $ U_{k}^{-2} $  in $ \alpha ( U_{i}^{*} U_{j} ) $ while ( \ref{lem3.2} ) and ( \ref{lem3.3} ) are obtained from the same coefficients in $ \alpha ( U_{i} U_{j}^{*} ) $ and $ \alpha ( U_{i} U_{j} )  $ respectively. \qed

      Now, by Lemma 2.12 in \cite{goswami} it follows that
            $ \tilde{\alpha}: C( \IT^{n} ) \otimes \clq \raro C( \IT^{n} ) \otimes \clq $ defined by                $\tilde{\alpha}(X \otimes Y)=\alpha(X)(1 \otimes Y)$ extends to  a unitary of the Hilbert $\clq$-module $L^2 (            C( \IT^{n} ) ,\tau ) \otimes \clq$ (or in other words, $\alpha$ extends to a unitary representation of $\clq$ on                 $L^2( C( \IT^{n} ) ,\tau)$). 
               But $\alpha$ keeps $W = {\rm Sp}\{U_{i},U_{i}^{*} : 1\leq i \leq n \}$ invariant.
               So $\alpha$ is a unitary representation of $\clq$ on $W$.Hence, the matrix ( say M ) corresponding to the $ 2n $ dimensional representation of $ \clq $ on $W$ is a unitary in $ M_{2n} ( \clq )$.
               
      From the definition of the action it follows that $ M =\left(  \begin {array} {cccc}
       A_{ij} &  B_{ij}^{*}  \\ B_{ij} & A_{ij}^{*} \end {array} \right )   $
      
      Since $ M $ is the matrix corresponding to a finite dimensional unitary representation, $ \kappa ( M_{k l } )= M^{-1}_{ k l } $ where $ \kappa $ denotes the antipode of $ \clq $ (See \cite{VanDaele})
      
    But $ M $ is a unitary, $ M^{-1} = M^{*} $
    
  So,$  ( k ( M_{k l} ) ) = \left ( \begin {array} {cccc}
   A_{ji}^{*} & B_{ji}^{*}  \\ B_{ji} & A_{ji} \end {array} \right ) $
      
  Now we apply the antipode $ \kappa $ to get some more relations.
  
  \blmma
    \label{Lemma 2}:
    
    $ \forall k $ and $ i \neq j $,
        \be A_{kj}^{*}A_{ki}^{*} = B_{kj} B_{ki} = A_{kj}^{*} B_{ki}^{*} = B_{kj} A_{ki} = B_{kj} A_{ki}^{*} = A_{kj} B_{ki} = 0  \ee
        \elmma
        
   {\it Proof :}\\ The result follows by applying $ \kappa $ on the equations $ A_{ik} A_{jk} = B_{ik} B_{jk} =  B^{*}_{ik} A_{jk} = A^{*}_{ik} B_{jk} = A_{ik}B_{jk} = B_{ik}A^{*}_{jk} = 0 $ obtained from Lemma  \ref{Lemma 3}.\qed
      
       \blmma
   \label{Lemma 5a}:

      $  A_{li} $ is a normal partial isometry $ \forall ~ l,i $  and hence has same domain and range.
      \elmma
      
   {\it Proof :}\\ From the relation ( \ref{lem1.1a} )  in Lemma \ref{Lemma 1}, we have by applying $ \kappa $, $ \sum (  A^{*}_{ji} A_{ji} +  B_{ji}B _{ji}^{*} ) = 1 $ .       
      Applying $ A_{li} $ on the right of this equation, we have 
      
    $  A_{li}^{*} A_{li} A_{li} + \sum_{j \neq l} ( A_{ji}^{*} A_{ji} A_{li} + B_{li} B_{li}^{*}A_{li} ) + \sum_{j \neq l} B_{ji} B_{ji}^{*}A_{li} = A_{li} .$     
   
   From Lemma \ref{Lemma 3}, we have $ A_{ji} A_{li} = 0 $ and $ B_{ji}^{*} A_{li} = 0 ~ \forall j \neq l $
   Moreover, from Lemma \ref{Lemma 1}, we have $ B_{li}^{*} A_{li} = 0 .$
   Applying these to the above equation, we have    
  \be \label {i}    A^{*}_{li} A_{li} A_{li} = A_{li}  \ee 
     
     Again, from the relation $  \sum_{j}( A_{ij} A^{*}_{ij} +  B_{ij}B^{*}_{ij} ) = 1 ~ \forall i $ in Lemma  \ref{Lemma 1}, applying $ \kappa $ and multiplying by $ A^{*}_{li} $ on the right, we have 
    $  A_{li} A_{li}^{*} A_{li}^{*} + \sum_{j \neq l} ( A_{ji} A^{*}_{ji} A^{*}_{li}  + B^{*}_{li} B_{li}A^{*}_{li}  + \sum_{j \neq l} B^{*}_{ji} B_{ji}A^{*}_{li} = A^{*}_{li} .$   
  From Lemma \ref{Lemma 3}, we have  $ A_{li} A_{ji} = 0 ~ \forall ~ j \neq l $( hence $  A^{*}_{ji} A^{*}_{li} = 0 $ ) and $ B_{ji} A^{*}_{li} = 0 .$ Moreover, we have $ B_{li} A^{*}_{li} = 0 $ from Lemma \ref{Lemma 1}.
  Hence, we have
     \be \label {ii}    A_{li} A^{*}_{li} A^{*}_{li} = A^{*}_{li}  \ee 
     
     From (\ref {i}), we have  \be \label {iii}  ( A^{*}_{li} A_{li} ) ( A_{li} A^{*}_{li} ) = A_{li} A^{*}_{li}  \ee 
     
     By taking $*$ on (\ref {ii}), we have \be \label {iv}  A_{li} A_{li} A^{*}_{li} = A_{li}  \ee
     
    Using this in (\ref {iii}), we have \be A_{li} A^{*}_{li} A_{li} = A_{li} A^{*}_{li} \ee
      and hence $ A_{li} $ is normal.
          
     So, $ A_{li} = A^{*}_{li} A_{li} A_{li} $ ( from ( \ref {i} ) )     
                $ = A_{li} A^{*}_{li} A_{li} $
                
      Therefore, $ A_{li} $ is a partial isometry which is normal and hence has same domain and range.\qed

      \blmma
   \label{Lemma 5b}:

   $ B_{li} $ is a normal partial isometry and hence has same domain and range.
   \elmma
   
   {\it Proof :}\\ 
   
 First we note that $ A_{ji} $ is a normal partial isometry and $ A_{ji} B_{li} = 0 ~ \forall ~ j \neq l $( obtained from Lemma \ref{Lemma 3} ) implies that  $ Ran ( A^{*}_{ji} ) \subseteq Ker ( B^{*}_{li} ) $ and hence $ Ran ( A_{ji} ) \subseteq Ker ( B^{*}_{li} ) $ which means  $ B^{*}_{li} A_{ji} = 0 ~ \forall j \neq l $.
   
   To obtain $ B^{*}_{li} B_{li} B_{li} = B_{li} $, we apply $ \kappa $ and multiply by $ B_{li} $ on the right of  ( \ref{lem1.5} ) and then use $ A^{*}_{li} B_{li} = 0 $ from Lemma \ref{Lemma 1} , $ A_{ji} B_{li} = 0 ~ \forall ~ j \neq l $( from Lemma \ref{Lemma 3}  which implies  $ B^{*}_{li} A_{ji} = 0 ~ \forall j \neq l $ as above ) and $ B_{ji} B_{li} = 0 ~ \forall j \neq l $ from Lemma \ref{Lemma 3} .
 
 Similarly, we have $ B_{li} B^{*}_{li} B^{*}_{li} = B^{*}_{li} $ by applying $ \kappa $ and multiplying by $ B^{*}_{li} $ on the right of ( \ref{lem1.1a} ) obtained from Lemma \ref{Lemma 1} and then use  $ A_{li} B^{*}_{li} = 0 $  ( Lemma \ref{Lemma 1} ), $ B_{li} A^{*}_{ji} = 0 ~ \forall ~ j \neq l .$ and  $ B_{li} B_{ji} = 0 ~ \forall j \neq l $ ( Lemma \ref{Lemma 3} ).

     Using   $ B^{*}_{li} B_{li} B_{li} = B_{li} $ and $ B_{li} B^{*}_{li} B^{*}_{li} = B^{*}_{li} $  as in Lemma \ref{Lemma 5a}, we have $ B_{li} $ is a normal partial isometry. \qed \vspace{4mm}

  Now, we use the condition  $ \alpha( U_{i} ) \alpha ( U_{j} ) = \alpha ( U_{j} ) \alpha ( U_{i} ) \forall i,j $
  
  \blmma
   \label{Lemma 6}:

      $ \forall k \neq l,$
  \be \label{lem6.1} A_{ik}A_{jl} + A_{il}A_{jk} = A_{jl}A_{ik} + A_{jk}A_{il} \ee
   
   \be A_{ik}B_{jl} + B_{il}A_{jk} = B_{jl}A_{ik} + A_{jk}B_{il} \ee
   
   \be B_{ik}A_{jl} + A_{il}B_{jk} = A_{jl}B_{ik} + B_{jk}A_{il} \ee
   
   \be B_{ik}B_{jl} + B_{il}B_{jk} = B_{jl}B_{ik} + B_{jk}B_{il}  \ee
   
  \elmma
  
 {\it Proof :}\\ The result follows by equating the coefficients of $ U_{k}U_{l}, U_{k}U^{-1}_{l}, U^{-1}_{k}U_{l} $ and $ U^{-1}_{k}U^{-1}_{l} $ ( where $  k \neq l $ ) in $ \alpha( U_{i} ) \alpha ( U_{j} ) = \alpha ( U_{j} ) \alpha ( U_{i} ) \forall i,j .$

    \qed
   
   \blmma
   \label{Lemma 7}:

   $ A_{ik}B_{jl} = B_{jl}A_{ik} \forall i \neq j, k \neq l $
   \elmma
   
   {\it Proof :}\\
               From Lemma \ref{Lemma 6},  we have $ \forall  k \neq l, A_{ik}B_{jl} - B_{jl}A_{ik} = A_{jk}B_{il} - B_{il}A_{jk} .$
  We consider the case where $ i \neq j.$ 
              
 We have,   $ Ran( A_{ik}B_{jl} - B_{jl}A_{ik} ) \subseteq Ran( A_{ik} ) + Ran(  B_{jl} ) \subseteq  Ran( B^{*}_{jl}B_{jl} + A^{*}_{ik}A_{ik} ) $ ( using the facts that $ A_{ik} $ and $ B_{jl} $ are normal partial isometries by Lemma \ref{Lemma 5a} and \ref{Lemma 5b} and also that  $ B^{*}_{jl}B_{jl} $ and $ A^{*}_{ik}A_{ik} $ are projections ).  
    
    Similarly, $ Ran ( A_{jk}B_{il} - B_{il}A_{jk} ) \subseteq  Ran( B^{*}_{il}B_{il} + A^{*}_{jk}A_{jk} ).$   
    
    Let \be T_{1} = A_{ik}B_{jl} - B_{jl}A_{ik} \ee
                   
         \be T_{2} = A_{jk}B_{il} - B_{il}A_{jk} \ee
          
        \be  T_{3} = B^{*}_{jl}B_{jl} + A^{*}_{ik}A_{ik} \ee
                   
         \be T_{4} =  B^{*}_{il}B_{il} + A^{*}_{jk}A_{jk} \ee
          
      Hence, $ T_{1} = T_{2}, Ran T_{1} \subseteq Ran T_{3}~, Ran T_{2}\subseteq Ran T_{4} .$             
                   
      We claim that $ T_{4}T_{3} = 0.$
      
      Then $ Ran ( T_{3} )  \subseteq Ker ( T_{4} ).$  
      
      But $ Ran T_{1} \subseteq Ran T_{3} $ will imply that $ Ran T_{1} \subseteq Ker T_{4}.$ Hence,$ Ran ( T_{2} ) \subseteq Ker ( T_{4} ) = { \overline{Ran ( T^{*}_{4})}}^{\bot} = { \overline{Ran ( T_{4})} }^{\bot} $
   But $ Ran ( T_{2} ) \subseteq Ran ( T_{4} ) $ which implies that $ Ran ( T_{2} ) = 0$ and hence both $ T_{2} $ and $ T_{1} $ are zero.
   Thus, the proof of the lemma will be complete if we can prove the claim.
   \be   T_{4} T_{3} = ( B^{*}_{il}B_{il} + A^{*}_{jk}A_{jk} ) ( B^{*}_{jl}B_{jl} + A^{*}_{ik}A_{ik} ) \ee
   
   \be = B^{*}_{il}B_{il}B^{*}_{jl}B_{jl} + B^{*}_{il}B_{il}A^{*}_{ik}A_{ik} + A^{*}_{jk}A_{jk}B^{*}_{jl}B_{jl} + A^{*}_{jk}A_{jk}A^{*}_{ik}A_{ik} \ee
   
   From Lemma \ref{Lemma 3}, we have $ \forall i \neq j , B_{il}B_{jl} = 0 $ implying $ B_{il} B^{*}_{jl} = 0 $ as $ B_{jl} $ is a normal partial isometry. 
  
  Again,from Lemma \ref{Lemma 2} $ \forall k \neq l, B_{il} A_{ik} = 0.$ Then $ A_{ik} $ is a normal partial isometry implies that $ B_{il} A^{*}_{ik} = 0 ~ \forall ~ k \neq l .$
  
  Similarly, by taking adjoint of the relation $ B_{jl} A^{*}_{jk} = 0 ~ \forall k \neq l $ obtained from Lemma \ref{Lemma 2}, we have $ A_{jk} B^{*}_{jl} = 0.$
  
  From Lemma \ref{Lemma 3}, we have $ A_{jk} A_{ik} = 0 ~ \forall ~ i \neq j .A_{ik} $ is a normal partial isometry implies that $ A_{jk} A^{*}_{ik} = 0 ~ \forall i \neq j .$
  
  Using these, we note that $ T_{4} T_{3} = 0 $ which proves the claim and hence the lemma.\qed

     \blmma
   \label{Lemma 8}:   
     
        \be A_{ik} B_{jk}= 0 =  B_{jk} A_{ik} \ee  
        
        \be A_{ki} B_{kj} = 0 = B_{kj} A_{ki} \ee  $ \forall i \neq j   $ and $ \forall k $
   \elmma

    {\it Proof :}\\           
  By Lemma \ref{Lemma 3}, we have $ A_{ik} B_{jk} = 0 $ and $ B_{jk} A^{*}_{ik} = 0 ~ \forall i \neq j $.  
  The second relation along with the fact that $ A_{ik} $ is a normal partial isometry implies that $ B_{jk} A{ik} = 0 \forall i \neq j .$
  
  Thus, $ A_{ik} B_{jk}= 0 =  B_{jk} A_{ik} ~ \forall i \neq j $
  
  Applying $ \kappa $ on the above equation and using $ B_{kj} $ and $ A_{ki} $ are normal partial isometries, we have $ A_{ki} B_{kj} = 0 = B_{kj} A_{ki} .$

   \qed \vspace{4mm}

        \blmma :
      \label{Lemma 9} 
     $ A_{ik} B_{ik} = B_{ik} A_{ik} ~ \forall i,k $
      \elmma 
      
     {\it Proof :}\\ 
     
     We have $ A^{*}_{ij} B_{ij} = 0 = B^{*}_{ij} A_{ij} $ from Lemma \ref{Lemma 1}.Using the fact that $ B_{ij} $ and $ A_{ij} $ are normal partial isometry we have $ A^{*}_{ij} B^{*}_{ij} = 0 = B^{*}_{ij} A^{*}_{ij} $ and hence $ A_{ij} B_{ij} = B_{ij} A_{ij} .$ \qed

    \blmma :
      \label{Lemma 10}

 $ A_{ik} A_{jl} = A_{jl} A_{ik} ~ \forall i \neq j, k \neq l .$
  
   \elmma

  {\it Proof :}\\ 
                               
 Using ( \ref{lem6.1} ) in Lemma \ref{Lemma 6},  we proceed as in Lemma \ref{Lemma 7} to get $ Ran( A_{ik}A_{jl} - A_{jl}A_{ik} ) \subseteq   Ran( A_{jl}A^{*}_{jl} + A_{ik}A^{*}_{ik} ) $ and  $ Ran ( A_{jk}A_{il} - A_{il}A_{jk} ) \subseteq  Ran( A_{il}A^{*}_{il} + A_{jk}A^{*}_{jk} ).$

  We claim that $ ( A_{ik}A^{*}_{ik} + A_{jl}A^{*}_{jl} )( A_{jk}A^{*}_{jk} + A_{il}A^{*}_{il} ) = 0.$
   
   Then by the same reasonings as given in Lemma \ref{Lemma 7} we will have : $  A_{jk} A_{il} = A_{il} A_{jk} .$

   To prove the claim,  we use $  A_{ik}A_{jk} = 0 ~ \forall ~ i \neq j $ from Lemma \ref{Lemma 3}( which implies $ A^{*}_{jk} A_{ik} = 0 ~ \forall ~ i \neq j  $ as $ A_{ik} $ is a normal partial isometry ), $ A^{*}_{il} A^{*}_{ik} = 0 ~ \forall ~ k \neq l $ from Lemma \ref{Lemma 2}( which implies $ A^{*}_{il} A_{ik} = 0 ~ \forall k \neq l $ as $ A_{ik} $ is a normal partial isometry ) and  $ A_{il} A_{jl} = 0 ~ \forall i \neq j .$ from Lemma \ref{Lemma 3} ( which implies  $ A^{*}_{jl} A_{il} = 0 ~ \forall i \neq j $ as $ A^{*}_{il} $ is a normal partial isometry ).

  \qed

    \blmma :
      \label{Lemma 11}

 \be A_{ik} A_{il} = A_{il} A_{ik} \forall  k \neq l \ee
 \be A_{ik} A_{jk} = A_{jk} A_{ik} \forall i \neq j \ee
   
   \elmma

  {\it Proof :}\\ 
  
        From Lemma \ref{Lemma 3}, we have $ A_{ki} A_{li} = 0 ~ \forall k \neq l .$
        
        Applying $ \kappa $ and taking adjoint, we have $ A_{ik} A_{il} = 0 ~ \forall k \neq l .$
        Interchanging $ k $ and $ l,$ we get $ A_{il}A_{ik} = 0 ~ \forall k \neq l .$
        Hence,$ A_{ik} A_{il} = A_{il} A_{ik} \forall  k \neq l .$ 
       
       From Lemma \ref{Lemma 3}, we have $ A_{ik} A_{jk} = 0 ~ \forall ~ i \neq j .$ Interchanging i and j, we have $ A_{jk} A_{ik} = 0 ~ \forall i \neq j .$ \qed
       
     \brmrk
      Proceeding in an exact similar way, we have that $ B_{ij} $ 's commute among themselves. 
                                                 \ermrk 
                                                 
     \bthm
     
      The Quantum Isometry group of $ \IT^{n} $ is commutative as a $ C^{\ast} $ algebra and hence coincides with the classical isometry group.
      
      \ethm
                                                  
    {\it Proof :}\\                                             
    Follows from the results in lemma \ref{Lemma 7} - \ref{Lemma 11} and the remark following them.  \qed

     \bcrlre

     Using Theorem 3.13 of \cite{bhowmick goswami}, we conclude that the Quantum Isometry Group of the noncommutative n tori $ \IT^{n}_{\theta} $       is a Rieffel deformation of the Quantum Isometry Group of $ \IT^{n} .$

     \ecrlre

  {\bf Acknowledgement :} We thank an anonymous referee for mentioning the paper \cite{VanDaele2} where A Van Daele removed  Woronowicz's separability assumption (in \cite{woro} ) for the $ C^* $ algebra of the underlying compact quantum group.

\end{document}